\documentclass[11pt]{amsart}
\usepackage{}
\usepackage{mathrsfs}
\usepackage{amsfonts}

\usepackage{amssymb,amsthm,amsmath,hyperref,txfonts}

\usepackage{graphicx,float}

\numberwithin{equation}{section}

\newtheorem{thm}{Theorem}[section]
\newtheorem{coro}{Corollary}[section]
\newtheorem{lem}{Lemma}[section]
\newtheorem{rem}{Remark}[section]
\newtheorem{prop}{Proposition}[section]

\theoremstyle{definition}
\newtheorem{defn}{Definition}[section]

\newcommand{\beq}{\begin{eqnarray}}
\newcommand{\eeq}{\end{eqnarray}}
\newcommand{\beqno}{\begin{eqnarray*}}
\newcommand{\eeqno}{\end{eqnarray*}}
\newcommand{\be}{\begin{equation}}
\newcommand{\ee}{\end{equation}}
\newcommand{\beno}{\begin{equation*}}
\newcommand{\eeno}{\end{equation*}}

\newcommand\tl{\tilde}

\newcommand\Dl{\Delta}
\newcommand\fr{\frac}

\newcommand\N{\mathbb{N}}
\newcommand\nb{\nabla}
\newcommand\nn{\nonumber}
\newcommand{\R}{\mathbb{R}}
\newcommand\Z{\mathbb{Z}}
\newcommand\sig{\sigma}

\newcommand\al{\alpha}
\newcommand{\dv}{\mathrm{div}}

\newcommand\lm{\lambda}

\newcommand\om{\omega}
\newcommand\pr{\partial}
\newcommand{\Rey}{\mbox{Re }}
\newcommand{\We}{\mbox{We }}
\newcommand{\ddl}{\dot{\Dl}_q}

\allowdisplaybreaks

\begin{document}
\title{Global solution to the incompressible Oldroyd-B model in hybrid Besov spaces}

\author{Ruizhao Zi}

\address{School of Mathematics and Statistics, Central
China Normal University, Wuhan 430079, China}
\email{ruizhao3805@163.com}



\begin{abstract}
This paper is dedicated to the Cauchy problem of the incompressible Oldroyd-B model with general coupling constant $\om\in (0,1)$.  It is shown that this set of equations admits a unique global solution in a certain hybrid Besov spaces for small initial data in $\dot{H}^s\cap\dot{B}^{\fr{d}{2}}_{2,1}$ with $-\fr{d}{2}<s<\fr{d}{2}-1$.  In particular, if $d\ge3,$ and taking $s=0$, then $\dot{H}^0\cap\dot{B}^{\fr{d}{2}}_{2,1}\approx B^{\fr{d}{2}}_{2,1}$. Since $B^{s}_{2,\infty}\hookrightarrow B^{\fr{d}{2}}_{2,1}, s>\fr{d}{2}$, this result extends the work by Chen and Miao [Nonlinear Anal.,{68}(2008), 1928--1939].
\end{abstract}

\maketitle
\section{Introduction}
We consider a typical model for viscoelastic fluids, the so called Oldroyd-B model \cite{Oldroyd58} in this paper. This type of fluids is  described by the following set of equations
\begin{eqnarray}\label{IOB}
\begin{cases}
\begin{array}{rrl}
u_t+(u\cdot\nabla) u-\eta_s\Delta u+\nabla
\Pi&=&\dv\tau,\\
\dv u&=&0,\\
\lambda(\tau_t+(u\cdot\nabla)\tau+g_\al(\tau,
\nabla u))+\tau&=&2\eta_e D(u),
\end{array}
\end{cases}
\end{eqnarray}
where $u$ and $\tau$ are the velocity and symmetric tensor of constrains of the fluids, respectively. $\Pi$ is the pressure which is the Lagrange multiplier for the divergence free condition. The quadratic form in $(\tau, \nb u)$ is given by $g_\al(\tau, \nabla u):=\tau W(u)-W(u)\tau-\al\left(D(u)\tau+\tau D(u)\right)$ for some $\al\in[-1,1]$, and
$D(u):=\frac{1}{2}(\nabla u+(\nabla u)^\top)$, $W(u):=\frac{1}{2}(\nabla u-(\nabla u)^\top)$ are the deformation tensor and the vorticity
tensor, respectively. Moreover, the parameter $\eta_s:=\eta\mu/\lm$ denotes the solvent viscosity, and $\eta_e:=\eta-\eta_s$ denotes the polymer viscosity, where $\eta$ is the total viscosity of the fluid, $\lm>0$ is the
relaxation time, and $\mu$ is the retardation time with $0<\mu<\lm$.

In the following, we would like to study system \eqref{IOB} in dimensionless variables, which takes the form
\begin{eqnarray}\label{IOBdimensionless}
\begin{cases}
\mathrm{Re}\left(u_t+(u\cdot\nabla) u\right)-(1-\om)\Delta u+\nabla
\Pi=\dv\tau,\\
\mathrm{We}(\tau_t+(u\cdot\nabla)\tau+g_\al(\tau,
\nabla u))+\tau=2\om D(u), \\
\dv u=0,
\end{cases}
\end{eqnarray}
with parameters Reynolds number $\Rey$,  Weissenberg number $\We$ and coupling constant $\om:=1-\fr{\mu}{\lm}\in(0,1)$ of the fluid .
For more details of the modeling, please refer to \cite{FZ14, GST10, Talhouk94} and references therein.

Some of the previous works in this direction can be summarized as follows. To our best knowledge, the incompressible Oldroyd-B model was firstly studied by Guillop\'e and Saut \cite{GS90}, where they obtained
a unique local strong solution to system \eqref{IOBdimensionless} in
suitable Sobolev spaces $H^s(\Omega)$ for the  situation of a {\em bounded domain} $\Omega \subset \R^3$. Moreover, this solution
is global provided both the data and the coupling
constant $\om$ between the two equations are sufficiently small.
The extensions to these results to the $L^p$-setting can be found in \cite{FGG98}. Similar results on {\em exterior domains} was established by Hieber, Naito and Shibata  \cite{HNS12}. The well-posedness  results in {\em scaling
invariant} Besov spaces on $\R^d, d\ge2$ were first given by Chemin and Masmoudi
\cite{CM01}.

All these results above were constructed under the assumption that the coupling constant $\om$ is small enough. This means that the coupling effect between the velocity $u$ and the symmetric tensor of constrains $\tau$ is weak and hence system \eqref{IOBdimensionless} corresponds closely to the classical incompressible Navier-Stokes equations.  From both the physical and mathematical point of view, it is more interesting to consider the strong coupling case, for which  the coupling constant $\om$ is not small. As a matter of fact,   the studies  in this direction have thrown up some interesting  results. For the situation of  {\em bounded domains},
the smallness restriction on the coupling constant $\om$ in
\cite{GS90} was removed by Molinet and Talhouk \cite{MT04}. As for the {\em exterior domains},   Fang, Hieber and the  author \cite{FHZ13} improved the main
result in \cite{HNS12} to the situation of {\em non-small}
coupling constant.   In the {\em whole space} $\mathbb{R}^d$ case,
Chen and Miao \cite{Chen-Miao08} obtained global solutions to  system \eqref{IOB} with small initial data in
$B^s_{2,\infty}, s>\frac{d}{2}$. For the critical $L^p$  framework, the smallness restriction on $\om$ in \cite{CM01} was removed by Fang, Zhang and the author \cite{Zi-Fang-Zhang14}  very recently. Existence of global {\em weak solutions} for large data and strong coupling was proved by Lions and Masmoudi in \cite{LM00} for the case $\al=0$. The general case $\al\neq0$ is still open up to now. For the Oldroyd-B fluids with {\em diffusive stress}, Constantin and Kliegl \cite{CK12} proved the global regularity of solutions in two dimensional case.

Besides, we would like to point out  that there are some other results on  Oldroyd-B fluids in the literature.
Indeed, Chemin and Masmoudi \cite{CM01} gave some {\em blow-up criterions} both for 2D and 3D cases. Later on, the 2D case was improved by Lei, Masmoudi and Zhou in \cite{LMZ10}. As for the 3D case, Kupferman, Mangoubi and Titi \cite{KMT08} established a Beale-Kato-Majda type blow-up criterion in terms of the $L^1_t(L^\infty_x)$ norm of $\tau$ in the {\em zero Reynolds number} regime. Further results,  describing the {\em incompressible limit problems} for Oldroyd-B fluids,  can be found in \cite{GST10,Lei06} for well-prepared initial data, and in \cite{FZ14} for ill-prepared initial data.    An approach based on {\em deformation tentor} was developed in \cite{HL14,Lei10,Lei07,LLZ08, LZ05,LLZ05,LZ08,Qian10, Zhang-Fang12}.

The aim of  this paper is to study the incompressible Oldroyd-B model \eqref{IOBdimensionless} with {\em non-small} coupling constant $\om$. We establish the  global solutions to system \eqref{IOBdimensionless} with small data $u_0$ and $\tau_0$ lying in $\dot{\mathcal{B}}^s\approx\dot{H}^s\cap\dot{B}^{\fr{d}{2}}_{2,1}, -\fr{d}{2}<s<\fr{d}{2}-1$. Like all the previous results \cite{Chen-Miao08, FHZ13, LM00, MT04} in $L^2$ framework with {\em non-small} coupling constant $\om$, the key point of the proof is to use the cancelation relation
\beno
(\dv\tau|u)+(D(u)|\tau)=0.
\eeno
The global estimates can be divided into two parts.  For the initial data in $\dot{B}^{\fr{d}{2}}_{2,1}$, owing to the Bernstein's inequality,  we can obtain both the smoothing effect of the velocity $u$ and the damping effect of the symmetric tensor of constrains $\tau$ in the high frequency case. While in the low frequency case,  the estimate fails to be true since $u$ and $\tau$ are treated as a whole, and $\|\ddl u\|_{L^2}+\|\ddl\tau\|_{L^2}$
can not be dominated by $\|\nb\ddl u\|_{L^2}+\|\ddl\tau\|_{L^2}$ any more (see \eqref{approx2} for details). In order to deal with the low frequency part, we impose an extra condition on the initial dada. This leads to the estimates for initial data in $\dot{H}^s$. It is worth noting that the estimates of nonlinear terms necessitate bounding the term $\|u\|_{L^2_t(\dot{B}^{\frac{d}{2}}_{2,1})}$. To do so, we decompose   $\|u\|_{L^2_t(\dot{B}^{\frac{d}{2}}_{2,1})}$ into $\|u\|^l_{L^2_t(\dot{B}^{\frac{d}{2}}_{2,1})}$ and $\|u\|^h_{L^2_t(\dot{B}^{\frac{d}{2}}_{2,1})}$. In particular, the low frequency part is bounded by $\|u\|^l_{L^2_t(\dot{H}^{s+1})}$; that is why we need $s<\fr{d}{2}-1$. Finally, combing the two parts estimates with initial data in $\dot{B}^{\fr{d}{2}}_{2,1}$ and $\dot{H}^s$, we obtain the global estimates for $(u, \tau)$.

\bigbreak\noindent{\bf Notations.}
For $s\in\R$, set
\beno
\|u\|^l_{\dot{B}^s_{2,1}}:=\sum_{q<0}2^{qs}\|\ddl u\|_{L^2},\quad \mathrm{and} \quad \|u\|^h_{\dot{B}^s_{2,1}}:=\sum_{q\ge0}2^{qs}\|\ddl u\|_{L^2}.
\eeno
Further more, let us denote by $\dot{B}^{s}_{\mathrm{h}}$ the space which consists of functions $u\in \mathcal{S}'$, such that $\|u\|^h_{\dot{B}^s_{2,1}}<\infty$. Throughout the paper, $C$ denotes various ``harmless'' positive constants, and
we sometimes use the notation $A \lesssim B$ as an equivalent to $A \le CB$. The
notation $A \approx B$ means that $A \lesssim B$ and $B \lesssim A$.

We shall obtain the existence and uniqueness of a solution $(u,\tau)$ to \eqref{IOBdimensionless} in the following space.

\begin{defn}\label{space}
For $T>0$, and $s\in\R$, let us denote
\beno
\mathcal{E}^s_T:=\left(\tl{C}_T(\dot{\mathcal{B}}^{s})\cap L^2_T(\dot{H}^{s+1})\cap L^1_T(\dot{B}^{\fr{d}{2}+1}_{\mathrm{h}})\right)^d\times\left(\tl{C}_T(\dot{\mathcal{B}}^{s})\cap L^2_T(\dot{H}^{s})\cap L^1_T(\dot{B}^{\fr{d}{2}}_{\mathrm{h}})\right)^{d\times d}.
\eeno
We use the notation $\mathcal{E}^s$ if $T=\infty$, changing $[0, T]$ into $[0,\infty)$ in the definition above. The definition of space $\dot{\mathcal{B}}^s$ can be found in Section \ref{tool}.
\end{defn}

Our main result reads as follows:
\begin{thm}\label{thm-global}
Let $d\ge2$, $-\fr{d}{2}<s<\fr{d}{2}-1$. Assume that $(u_0,\tau_0)\in \left(\dot{\mathcal{B}}^{s}\right)^d\times \left(\dot{\mathcal{B}}^{s}\right)^{d\times d}$ with $\dv u_0=0$.
There exist two positive constants $c$ and $M$, depending on $s, d, \om, \Rey$ and $\We$, such that if
                                      \beno
                                      \|u_0\|_{\dot{\mathcal{B}}^{s}}+\|\tau_0\|_{\dot{\mathcal{B}}^{s}}\le c,
                                      \eeno
system \eqref{IOBdimensionless} admits a unique global solution $(u,\tau)$ in $\mathcal {E}^{s}$ with
                                    \beno
                                     \|(u,\tau)\|_{\mathcal{E}^{s}}\leq M\left(\|u_0\|_{\dot{\mathcal{B}}^{s}}+\|\tau_0\|_{\dot{\mathcal{B}}^{s}}\right).
                                    \eeno
\end{thm}

\section{The Functional Tool Box}\label{tool}
\noindent The results of the present paper rely on the use of  a
dyadic partition of unity with respect to the Fourier variables, the so-called the
\textit{Littlewood-Paley  decomposition}. Let us briefly explain how
it may be built in the case $x\in \R^d$ which the readers may see more details
in \cite{Bahouri-Chemin-Danchin11,Ch1}. Let $(\chi, \varphi)$ be a couple of $C^\infty$ functions satisfying
$$\hbox{Supp}\chi\subset\{|\xi|\leq\frac{4}{3}\},
\ \ \ \
\hbox{Supp}\varphi\subset\{\frac{3}{4}\leq|\xi|\leq\frac{8}{3}\},
$$
and
$$\chi(\xi)+\sum_{q\geq0}\varphi(2^{-q}\xi)=1,$$

$$\sum_{q\in \mathbb{Z}}\varphi(2^{-q}\xi)=1, \quad \textrm{for} \quad \xi\neq0.$$
Set $\varphi_q(\xi)=\varphi(2^{-q}\xi),$
$h_q=\mathcal{F}^{-1}(\varphi_q),$ and
$\tilde{h}=\mathcal{F}^{-1}(\chi)$. The dyadic blocks and the low-frequency cutoff operators are defined for all $q\in\mathbb{Z}$ by
$$\dot{\Delta}_{q}u=\varphi(2^{-q}\mathrm{D})u=\int_{\R^d}h_q(y)u(x-y)dy,$$
$$\dot{S}_qu=\chi(2^{-q}\mathrm{D})u=\int_{\R^d}\tl{h}_q(y)u(x-y)dy.$$
Then
\begin{equation}\label{e2.1}
u=\sum_{q\in \mathbb{Z}}\Delta_qu,
\end{equation}
holds for tempered distributions {\em modulo polynomials}. As working modulo polynomials is not appropriate for nonlinear problems, we
shall restrict our attention to the set $\mathcal {S}'_h$ of tempered distributions $u$ such that
$$
\lim_{q\rightarrow-\infty}\|\dot{S}_qu\|_{L^\infty}=0.
$$
Note that \eqref{e2.1} holds true whenever $u$ is in $\mathcal{S}'_h$ and that one may write
$$
\dot{S}_qu=\sum_{p\leq q-1}\dot{\Dl}_{p}u.
$$
Besides, we would like to mention that the Littlewood-Paley decomposition
has a nice property of quasi-orthogonality:
\begin{equation}\label{e2.2}
\dot{\Delta}_p\dot{\Delta}_qu\equiv 0\ \ \hbox{if}\ \ \ |p-q|\geq 2\ \
\hbox{and}\ \ \dot{\Delta}_p(\dot{S}_{q-1}u\dot{\Delta}_qu)\equiv 0\ \ \hbox{if}\ \ \
|p-q|\geq 5.
\end{equation}
One can now  give the definition of
homogeneous Besov spaces.
\begin{defn}\label{D2.1}
For $s\in\R$, $(p,r)\in[1,\infty]^2$, and
$u\in\mathcal{S}'(\R^d),$ we set
$$\|u\|_{\dot{B}_{p,r}^s}=\left\|2^{ sq}\|\dot{\Delta}_qu\|_{L^p} \right\|_{\ell^r}.$$
We then define the space
$\dot{B}_{p,r}^s:=\{u\in\mathcal{S}'_h(\R^d),\
\|u\|_{\dot{B}_{p,r}^s}<\infty\}$.
\end{defn}
\begin{rem}

The inhomogeneous Besov spaces can be defined in a similar way. Indeed, for $u\in \mathcal {S}'(\R^d)$,
we set
$$\Delta_qu=0\ \ \mathrm{if}\ \ q<-1,\ \ \ \Delta_{-1}u=\chi(\mathrm{D})u,$$
$$\Delta_{q}u=\varphi(2^{-q}\mathrm{D})u\ \ \mathrm{if}\ \ q\geq0, \quad \mathrm{and} \quad S_qu=\sum_{p\leq q-1}\Delta_pu.$$
Then for all $u\in \mathcal {S}'(\R^d)$, we have the inhomogeneous
Littlewood-Paley decomposition $u = \sum_{q\in\mathbb{Z}} \Dl_q u$,
and for $(p, r ) \in [1,+\infty]^2, s \in \R$, we define the
inhomogeneous Besov space $B_{p,r}^s$ as
$$B_{p,r}^s=\{u\in\mathcal{S}'(\R^d), \|u\|_{B_{p,r}^s}:=\left\|2^{ sq}\|\Delta_qu\|_{L^p}\right\|_{\ell^r}<\infty\}$$
\end{rem}
We also need the following hybrid Besov space in this paper.

\begin{defn}\label{D2.2}
For $s\in\R$, and
$u\in\mathcal{S}'(\R^d),$ we set
$$\|u\|_{\dot{\mathcal{B}}^s}=\left(\sum_{q<0}2^{2qs}\|\ddl u\|_{L^2}^2\right)^{\fr12}+\sum_{q\ge0}2^{q\fr{d}{2}}\|\ddl u\|_{L^2}.$$
We then define the space
$\dot{\mathcal{B}}^s:=\{u\in\mathcal{S}'_h(\R^d),\
\|u\|_{\dot{\mathcal{B}}^s}<\infty\}$.
\end{defn}

\begin{rem}
It is easy to verify that
\beno
\dot{\mathcal{B}}^s\approx \dot{H}^s\cap\dot{B}^{\fr{d}{2}}_{2,1},
\eeno
provided $s<\fr{d}{2}$.
\end{rem}

The following lemma describes the way derivatives act on spectrally localized functions.
\begin{lem}[Bernstein's inequalities]\label{Bernstein}
Let $k\in\N$ and $0<r<R$. There exists a constant $C$ depending on $r, R$ and $d$ such that for all $(a,b)\in[1,\infty]^2$, we have for all $\lm>0$ and multi-index $\al$
\begin{itemize}
\item If $\mathrm{Supp} \hat{f}\subset B(0,\lm R)$, then $\sup_{\al=k}\|\pr^\al f\|_{L^b}\le C^{k+1}\lm^{k+d(\fr1a-\fr1b)}\|f\|_{L^a}$.
\item If $\mathrm{Supp} \hat{f}\subset \mathfrak{C}(0,\lm r, \lm R)$, then $C^{-k-1}\lm^k\|f\|_{L^a}\le\sup_{|\al|=k}\|\pr^\al f\|_{L^a}\le C^{k+1}\lm^k\|f\|_{L^a}$
\end{itemize}
\end{lem}

Next we  recall a few nonlinear estimates in Besov spaces which may be
obtained by means of paradifferential calculus. Firstly introduced
 by J. M. Bony in \cite{Bony81}, the paraproduct between $f$
and $g$ is defined by
$$\dot{T}_fg=\sum_{q\in\mathbb{Z}}\dot{S}_{q-1}f\dot{\Delta}_qg,$$
and the remainder is given by
$$\dot{R}(f,g)=\sum_{q\geq -1}\dot{\Delta}_qf\tilde{\dot{\Delta}}_qg$$
with
$$\tilde{\dot{\Delta}}_qg=(\dot{\Delta}_{q-1}+\dot{\Delta}_{q}+\dot{\Delta}_{q+1})g.$$
We have the following so-called Bony's decomposition:
 \be\label{Bony-decom}
fg=\dot{T}_fg+\dot{T}_gf+\dot{R}(f,g).
 \ee
The paraproduct $\dot{T}$ and the remainder $\dot{R}$ operators satisfy the following
continuous properties.

\begin{prop}\label{p-TR}
For all $s\in\mathbb{R}$, $\sigma>0$, and $1\leq p, p_1, p_2, r, r_1, r_2\leq\infty,$ the
paraproduct $\dot T$ is a bilinear, continuous operator from $L^{\infty}\times \dot{B}_{p,r}^s$ to $
\dot{B}_{p,r}^{s}$ and from $\dot{B}_{\infty,r_1}^{-\sigma}\times \dot{B}_{p,r_2}^s$ to
$\dot{B}_{p,r}^{s-\sigma}$  with $\frac{1}{r}=\min\{1, \frac{1}{r_1}+\frac{1}{r_2}\}$. The remainder $\dot R$ is bilinear continuous from
$\dot{B}_{p_1,r_1}^{s_1}\times \dot{B}_{p_2,r_2}^{s_2}$ to $
\dot{B}_{p,r}^{s_1+s_2}$ with
$s_1+s_2>0$,  $\frac{1}{p}=\frac{1}{p_1}+\frac{1}{p_2}\leq1$, and $\frac{1}{r}=\frac{1}{r_1}+\frac{1}{r_2}\leq1$.
\end{prop}
In view of \eqref{Bony-decom}, Proposition \ref{p-TR} and Bernstein's inequalities,  one easily deduces the following  product estimates:
\begin{coro}\label{coro-product}
There hold:
\be\label{product1}
\|uv\|_{\dot{H}^{s}}\leq C\|u\|_{\dot{B}^{\fr{d}{2}}_{2,1}}\|v\|_{\dot{H}^{s}},\quad\mathrm{if}\quad s\in(-\fr{d}{2},\fr{d}{2}). 
\ee
\be\label{product1.5}
\|uv\|_{\dot{H}^{s}}\leq C\|u\|_{\dot{H}^{s+1}}\|v\|_{\dot{B}^{\fr{d}{2}-1}_{2,\infty}},\quad\mathrm{if}\quad s\in(-\fr{d}{2},\fr{d}{2}-1).
\ee
and
\be\label{product2}
\|uv\|_{\dot{B}^{\fr{d}{2}}_{2,1}}\leq C\|u\|_{\dot{B}^{\fr{d}{2}}_{2,1}}\|v\|_{\dot{B}^{\fr{d}{2}}_{2,1}}.
\ee
\end{coro}
The study of non-stationary PDEs requires spaces of the type
$L^\rho_T(X)=L^\rho(0,T;X)$ for appropriate Banach spaces $X$. In
our case, we expect $X$ to be a  Besov space, so that it
is natural to localize the equations through Littlewood-Paley
decomposition. We then get estimates for each dyadic block and
perform integration in time. But, in doing so, we obtain the bounds
in spaces which are not of the type $L^\rho(0,T;\dot{B}^s_{p,r})$. That
 naturally leads to the following definition introduced by Chemin and Lerner in \cite{CL}.
\begin{defn}\label{defn-chemin-lerne}
For $\rho\in[1,+\infty]$, $s\in\R$, and $T\in(0,+\infty)$, we set
$$\|u\|_{\tilde{L}^\rho_T(\dot{B}^s_{p,r})}=\left\|2^{qs}
\|\dot{\Delta}_qu(t)\|_{L^\rho_T(L^p)}\right\|_{\ell^r}
$$
and denote by
$\tilde{L}^\rho_T(\dot{B}^s_{p,r})$ the subset of distributions
$u\in\mathcal{D}'([0,T]; \mathcal{S}'_h (\mathbb{R}^N))$ with finite
$\|u\|_{\tilde{L}^\rho_T(\dot{B}^s_{p,r})}$ norm. When $T=+\infty$, the index $T$ is
omitted. We
further denote $\tilde{C}_T(\dot{B}^s_{p,r})=C([0,T];\dot{B}^s_{p,r})\cap
\tilde{L}^\infty_{T}(\dot{B}^s_{p,r}) $.
 \end{defn}
\begin{rem}\label{rem-CL-holder}
All the properties of continuity for the paraproduct, remainder, and product remain true for the Chemin-Lerner spaces. The exponent $\rho$ just has to behave according to H\"{o}lder's ineauality for the time variable.
\end{rem}

\begin{rem}\label{rem-CL-minkowski}
The spaces $\tl{L}^\rho_T(\dot{B}^s_{p,r})$ can be linked with the classical space $L^\rho_T(\dot{B}^s_{p,r})$ via the Minkowski inequality:
\beno
\|u\|_{\tl{L}^\rho_T(\dot{B}^s_{p,r})}\le\|u\|_{L^\rho_T(\dot{B}^s_{p,r})}\quad \mathrm{if}\quad r\ge\rho,\qquad \|u\|_{\tl{L}^\rho_T(\dot{B}^s_{p,r})}\ge\|u\|_{L^\rho_T(\dot{B}^s_{p,r})}\quad \mathrm{if}\quad r\le\rho.
\eeno
\end{rem}

\section{Global existence}\label{Sec-E}
\noindent In order to construct the global solutions to the incompressible Oldroyd-B model \eqref{IOBdimensionless},  we shall used the classical
Friedrichs method to approximate the system \eqref{IOBdimensionless} by a cut-off in the frequency space. Noting that this method has been applied to Oldroyd-B model in \cite{CM01,Chen-Miao08} before, to avoid unnecessary repetition,  we omit the details of  approximation  in this paper. In the following,  the global estimates of $(u, \tau)$ will be given directly. To begin with, let us first of all localize the system \eqref{IOBdimensionless} as follows,
 \beq\label{loc_fre}
\begin{cases}
2\om\mathrm{Re}\dot{\Dl}_qu_t-2\om (1-\om)\dot{\Dl}_q\Dl u+2\om\nb\dot{\Dl}_q\Pi=2\om\dv \dot{\Dl}_q\tau-2\om\mathrm{Re}\dot{\Dl}_q(u\cdot\nb u),\\
\mathrm{We}\left(\dot{\Dl}_q\tau_t+u\cdot\nb\dot{\Dl}_q\tau\right)+\dot{\Dl}_q\tau=2\om D(\dot{\Dl}_qu)-\mathrm{We}\left([\dot{\Dl}_q,u]\cdot\nb\tau+\dot{\Dl}_qg_\al(\tau,\nb u)\right).
\end{cases}
 \eeq
Taking the $L^2$ inner product of $\eqref{loc_fre}_1$ and $\eqref{loc_fre}_2$ with $\ddl u$ and $\ddl\tau$, respectively, using the relation $(\dv \tau|u)+(D(u)|\tau)=0$ and the divergence free condition of $u$, we obtain
 \beq\label{est_loc}
\nn&&\frac12\frac{d}{dt}\left(2\om\mathrm{Re}\|\ddl u\|_{L^2}^2+\mathrm{We}\|\ddl\tau\|_{L^2}^2\right)+2\om(1-\om)\|\nb\ddl u\|_{L^2}^2+\|\ddl\tau\|_{L^2}^2\\
&\le&2\om\mathrm{Re}\|\ddl(u\cdot\nb u)\|_{L^2}\|\ddl u\|_{L^2}+\mathrm{We}\left(\|[\ddl,u]\cdot\nb\tau\|_{L^2}+\|\ddl g_\al(\tau,\nb u)\|_{L^2}\right)\|\ddl \tau\|_{L^2}.
 \eeq
It follows that
 \beq\label{H^s_est}
\nn&&\om\mathrm{Re}\|u\|_{\tl{L}^\infty_t(\dot{H}^s)}^2+\frac{\mathrm{We}}{2}\|\tau\|_{\tl{L}^\infty_t(\dot{H}^s)}^2+2\om(1-\om)\|\nb u\|_{L^2_t(\dot{H}^s)}^2+\|\tau\|_{L^2_t(\dot{H}^s)}^2\\
\nn&\le&\om\mathrm{Re}\|u_0\|_{\dot{H}^s}^2+\frac{\mathrm{We}}{2}\|\tau_0\|_{\dot{H}^s}^2+2\om\mathrm{Re}\int_0^t\sum_{q\in\mathbb{Z}}
2^{2qs}\|\ddl(u\cdot\nb u)\|_{L^2}\|\ddl u\|_{L^2}dt'\\
&&+\mathrm{We}\int_0^t\sum_{q\in\mathbb{Z}}2^{2qs}\left(\|[\ddl, u]\cdot\nb\tau\|_{L^2}+\|\ddl g_\al(\tau,\nb u)\|_{L^2}\right)\|\ddl\tau\|_{L^2}dt'.
 \eeq
Now we estimate the last three terms on the righthand side of \eqref{H^s_est} one by one. In fact, in view of H\"{o}der's inequality and the product estimate \eqref{product1}, we infer that for $-\frac{d}{2}<s<\frac{d}{2}$, there holds
 \beq\label{est_convu}
\nn&&2\om\mathrm{Re}\int_0^t\sum_{q\in\mathbb{Z}}
2^{2qs}\|\ddl(u\cdot\nb u)\|_{L^2}\|\ddl u\|_{L^2}dt'\\
\nn&\leq&2\om\mathrm{Re}\int_0^t\|u\cdot\nb u\|_{\dot{H}^s}\|u\|_{\dot{H}^s}dt'\\
\nn&\leq&C\om\mathrm{Re}\int_0^t\|\nb u\|_{\dot{H}^s}\|u\|_{\dot{B}^{\frac{d}{2}}_{2,1}}\|u\|_{\dot{H}^s}dt'\\
&\leq&C\om\mathrm{Re}\|u\|_{L^\infty_t(\dot{H}^s)}\|\nb u\|_{L^2_t(\dot{H}^s)}\|u\|_{L^2_t(\dot{B}^{\frac{d}{2}}_{2,1})}.
 \eeq
Noting that if $s<\frac{d}{2}-1$, we can bound $\|u\|_{L^2_t(\dot{B}^{\frac{d}{2}}_{2,1})}$ as follows:
 \beq\label{est_u}
\|u\|_{L^2_t(\dot{B}^{\frac{d}{2}}_{2,1})}\nn&\le&\|u\|^l_{\tl{L}^2_t(\dot{B}^{\frac{d}{2}}_{2,1})}+\|u\|^h_{\tl{L}^2_t(\dot{B}^{\frac{d}{2}}_{2,1})}\\
\nn&\le&C\left(\|u\|_{L^2_t(\dot{H}^{s+1})}^l+\|u\|^h_{\tl{L}^2_t(\dot{B}^{\frac{d}{2}+\frac12}_{2,1})}\right)\\
&\le&C\left(\|u\|_{L^2_t(\dot{H}^{s+1})}+\left(\|u\|^h_{L^\infty_t(\dot{B}^{\frac{d}{2}}_{2,1})}
\|u\|^h_{L^1_t(\dot{B}^{\frac{d}{2}+1}_{2,1})}\right)^{\frac12}\right)
 \eeq
Inserting \eqref{est_u} into \eqref{est_convu}, we arrive at
  \beq\label{est_conv-u}
\nn&&2\om\mathrm{Re}\int_0^t\sum_{q\in\mathbb{Z}}
2^{2qs}\|\ddl(u\cdot\nb u)\|_{L^2}\|\ddl u\|_{L^2}dt'\\
&\leq&C\om\mathrm{Re}\|u\|_{L^\infty_t(\dot{H}^s)}\|\nb u\|_{L^2_t(\dot{H}^s)}\left(\|u\|_{L^2_t(\dot{H}^{s+1})}+\left(\|u\|^h_{L^\infty_t(\dot{B}^{\frac{d}{2}}_{2,1})}
\|u\|^h_{L^1_t(\dot{B}^{\frac{d}{2}+1}_{2,1})}\right)^{\frac12}\right),
 \eeq
with $-\frac{d}{2}<s<\frac{d}{2}-1$. Similarly, we have
 \beq\label{est_g}
\nn&&\mathrm{We}\int_0^t\sum_{q\in\mathbb{Z}}2^{2qs}\|\ddl g_\al(\tau,\nb u)\|_{L^2}\|\ddl\tau\|_{L^2}dt'\\
\nn&\leq&\mathrm{We}\int_0^t\| g_\al(\tau,\nb u)\|_{\dot{H}^s}\|\tau\|_{\dot{H}^s}dt'\\
\nn&\leq&C\mathrm{We}\int_0^t\|\tau\|_{\dot{B}^{\frac{d}{2}}_{2,1}}\|\nb u\|_{\dot{H}^s}\|\tau\|_{\dot{H}^s}dt'\\
&\leq&C\mathrm{We}\|\tau\|_{\tl{L}^\infty_t(\dot{B}^{\frac{d}{2}}_{2,1})}\|\nb u\|_{L^2_t(\dot{H}^s)}\|\tau\|_{L^2_t(\dot{H}^s)}, \quad\hbox{for}\quad  -\frac{d}{2}<s<\frac{d}{2}.
 \eeq
Finally, using H\"{o}lder's inequality and the commutator estimate, c. f. \cite{Bahouri-Chemin-Danchin11}, for $-\frac{d}{2}-1<s<\frac{d}{2}$, we are led to
 \beq\label{est_com-tau}
\nn&&\mathrm{We}\int_0^t\sum_{q\in\mathbb{Z}}2^{2qs}\|[\ddl, u]\cdot\nb\tau\|_{L^2}\|\ddl\tau\|_{L^2}dt'\\
\nn&\leq&\mathrm{We}\int_0^t\left(\sum_{q\in\mathbb{Z}}2^{2qs}\|[\ddl,u]\cdot\nb\tau\|^2_{L^2}\right)^{\frac12}\|\tau\|_{\dot{H}^s}dt'\\
\nn&\leq&C\mathrm{We}\int_0^t\|\nb u\|_{\dot{B}^{\frac{d}{2}}_{2,1}}\|\tau\|_{\dot{H}^s}^2dt'\\
\nn&=&C\mathrm{We}\int_0^t\|\nb u\|^l_{\dot{B}^{\frac{d}{2}}_{2,1}}\|\tau\|_{\dot{H}^s}^2+\|\nb u\|^h_{\dot{B}^{\frac{d}{2}}_{2,1}}\|\tau\|_{\dot{H}^s}^2dt'\\
\nn&\leq&C\mathrm{We}\|\nb u\|^l_{L^2_t(\dot{B}^{\frac{d}{2}}_{2,1})}\|\tau\|_{L^\infty_t(\dot{H}^s)}\|\tau\|_{L^2_t(\dot{H}^s)}
                   +C\mathrm{We}\|\nb u\|^h_{L^1_t(\dot{B}^{\frac{d}{2}}_{2,1})}\|\tau\|_{L^\infty_t(\dot{H}^s)}^2\\
&\leq&C\mathrm{We}\|\nb u\|_{L^2_t(\dot{H}^{s})}\|\tau\|_{L^\infty_t(\dot{H}^s)}\|\tau\|_{L^2_t(\dot{H}^s)}
                   +C\mathrm{We}\|\nb u\|^h_{L^1_t(\dot{B}^{\frac{d}{2}}_{2,1})}\|\tau\|_{L^\infty_t(\dot{H}^s)}^2.
 \eeq
Substituting \eqref{est_conv-u}--\eqref{est_com-tau} into \eqref{H^s_est} yields
 \beq\label{est_H^s}
\nn&&\om\mathrm{Re}\|u\|_{\tl{L}^\infty_t(\dot{H}^s)}^2+\frac{\mathrm{We}}{2}\|\tau\|_{\tl{L}^\infty_t(\dot{H}^s)}^2+2\om(1-\om)\|\nb u\|_{L^2_t(\dot{H}^s)}^2+\|\tau\|_{L^2_t(\dot{H}^s)}^2\\
\nn&\le&\om\mathrm{Re}\|u_0\|_{\dot{H}^s}^2+\frac{\mathrm{We}}{2}\|\tau_0\|_{\dot{H}^s}^2+C\mathrm{We}\|\tau\|_{\tl{L}^\infty_t(\dot{B}^{\frac{d}{2}}_{2,1})}\|\nb u\|_{L^2_t(\dot{H}^s)}\|\tau\|_{L^2_t(\dot{H}^s)}\\
\nn&&+C\mathrm{We}\|\nb u\|_{L^2_t(\dot{H}^{s})}\|\tau\|_{L^\infty_t(\dot{H}^s)}\|\tau\|_{L^2_t(\dot{H}^s)}
                   +C\mathrm{We}\|\nb u\|^h_{L^1_t(\dot{B}^{\frac{d}{2}}_{2,1})}\|\tau\|_{L^\infty_t(\dot{H}^s)}^2\\
&&+C\om\mathrm{Re}\|u\|_{L^\infty_t(\dot{H}^s)}\|\nb u\|_{L^2_t(\dot{H}^s)}\left(\|u\|_{L^2_t(\dot{H}^{s+1})}+\left(\|u\|^h_{L^\infty_t(\dot{B}^{\frac{d}{2}}_{2,1})}
\|u\|^h_{L^1_t(\dot{B}^{\frac{d}{2}+1}_{2,1})}\right)^{\frac12}\right).
 \eeq
To close \eqref{est_H^s}, we have to estimate the high frequency part of $u$ and $\tau$ in the space $\dot{B}^{\frac{d}{2}}_{2,1}$. To this end, we first notice that
 \be\label{approx1}
2\om\mathrm{Re}\|\ddl u\|_{L^2}^2+\mathrm{We}\|\ddl\tau\|_{L^2}^2\approx\left(\sqrt{2\om\mathrm{Re}}\|\ddl u\|_{L^2}+\sqrt{\mathrm{We}}\|\ddl\tau\|_{L^2}\right)^2,
 \ee
and for $q\ge0$,
 \beq\label{approx2}
\nn2\om(1-\om)\|\nb\ddl u\|_{L^2}^2+\|\ddl\tau\|_{L^2}^2&\approx&2\om(1-\om)2^{2q}\|\ddl u\|_{L^2}^2+\|\ddl\tau\|_{L^2}^2\\
\nn&\approx&\left(\sqrt{2\om(1-\om)}2^q\|\ddl u\|_{L^2}+\|\ddl\tau\|_{L^2}\right)^2\\
\nn&\ge&\min\left\{\sqrt{\frac{1-\om}{\mathrm{Re}}},\frac{1}{\sqrt{\mathrm{We}}}\right\}\left(\sqrt{2\om(1-\om)}2^q\|\ddl u\|_{L^2}+\|\ddl\tau\|_{L^2}\right)\\
&&\times\left(\sqrt{2\om\mathrm{Re}}\|\ddl u\|_{L^2}+\sqrt{\mathrm{We}}\|\ddl\tau\|_{L^2}\right).
 \eeq
It follows from \eqref{est_loc}, \eqref{approx1} and  \eqref{approx2} that, if $q\ge0$, there holds
 \beq
\nn&&\frac{d}{dt}\sqrt{2\om\mathrm{Re}\|\ddl u\|_{L^2}^2+\mathrm{We}\|\ddl\tau\|_{L^2}^2}\\
\nn&&\qquad\qquad\qquad\qquad+\min\left\{\sqrt{\frac{1-\om}{\mathrm{Re}}},\frac{1}{\sqrt{\mathrm{We}}}\right\}\left(\sqrt{2\om(1-\om)}2^q\|\ddl u\|_{L^2}+\|\ddl\tau\|_{L^2}\right)\\
&\lesssim&\sqrt{2\om\mathrm{Re}}\|\ddl(u\cdot\nb u)\|_{L^2}+\sqrt{\mathrm{We}}\left(\|[\ddl, u]\cdot\nb\tau\|_{L^2} +\|\ddl g_\al(\tau,\nb u)\|_{L^2}\right).
 \eeq
Integrating the above equation with respect to $t$, multiplying the resulting inequality by $2^{q{\frac{d}{2}}}$, and then summing w. r. t. $q$ over all the nonnegative integers, we find that
 \beq\label{B_est}
\nn&&\sqrt{2\om\mathrm{Re}}\|u\|^h_{\tl{L}^\infty_t(\dot{B}^{\frac{d}{2}}_{2,1})}+\sqrt{\mathrm{We}}\|\tau\|^h_{\tl{L}^\infty_t(\dot{B}^{\frac{d}{2}}_{2,1})}\\
\nn&&\qquad\qquad\qquad+\min\left\{\sqrt{\frac{1-\om}{\mathrm{Re}}},\frac{1}{\sqrt{\mathrm{We}}}\right\}\left(\sqrt{2\om(1-\om)}\| u\|^h_{{L}^1_t(\dot{B}^{\frac{d}{2}+1}_{2,1})}+\|\tau\|^h_{{L}^1_t(\dot{B}^{\frac{d}{2}}_{2,1})}\right)\\
\nn&\lesssim&\sqrt{2\om\mathrm{Re}}\|u_0\|^h_{\dot{B}^{\frac{d}{2}}_{2,1}}+\sqrt{\mathrm{We}}\|\tau_0\|^h_{\dot{B}^{\frac{d}{2}}_{2,1}}
+\sqrt{2\om\mathrm{Re}}\int_0^t\|u\cdot\nb u\|_{\dot{B}^{\frac{d}{2}}_{2,1}}dt'\\
&&+\sqrt{\mathrm{We}}\int_0^t\sum_{q\in\mathbb{Z}}2^{q\frac{d}{2}}\|[\ddl,u]\cdot\nb\tau\|_{L^2}+\|g_\al(\tau,\nb u)\|_{\dot{B}^{\frac{d}{2}}_{2,1}}dt'.
 \eeq
Product estimate \eqref{product2} in Besov space  and \eqref{est_u} imply that, for $-\frac{d}{2}<s<\frac{d}{2}-1$,
 \beq\label{est_conv-u-B}
\nn&&\sqrt{2\om\mathrm{Re}}\int_0^t\|u\cdot\nb u\|_{\dot{B}^{\frac{d}{2}}_{2,1}}dt'\\
\nn&\leq&C\sqrt{2\om\mathrm{Re}}\int_0^t\|u\|_{\dot{B}^{\frac{d}{2}}_{2,1}}\|\nb u\|_{\dot{B}^{\frac{d}{2}}_{2,1}}dt'\\
\nn&=&C\sqrt{2\om\mathrm{Re}}\int_0^t\|u\|_{\dot{B}^{\frac{d}{2}}_{2,1}}\left(\|\nb u\|^l_{\dot{B}^{\frac{d}{2}}_{2,1}}+\|\nb u\|^h_{\dot{B}^{\frac{d}{2}}_{2,1}}\right)dt'\\
\nn&\le&C\sqrt{2\om\mathrm{Re}}\left(\|u\|_{L^\infty_t(\dot{B}^{\frac{d}{2}}_{2,1})}\|\nb u\|^h_{L^1_t(\dot{B}^{\frac{d}{2}}_{2,1})}+\|u\|_{L^2_t(\dot{B}^{\frac{d}{2}}_{2,1})}\|\nb u\|_{L^2_t(\dot{B}^{\frac{d}{2}}_{2,1})}^l\right)\\
&\le&C\sqrt{2\om\mathrm{Re}}\left(\|u\|_{L^\infty_t(\dot{B}^{\frac{d}{2}}_{2,1})}\|\nb u\|^h_{L^1_t(\dot{B}^{\frac{d}{2}}_{2,1})}+\left(\|u\|_{L^2_t(\dot{H}^{s+1})}+\left(\|u\|^h_{L^\infty_t(\dot{B}^{\frac{d}{2}}_{2,1})}
\|u\|^h_{L^1_t(\dot{B}^{\frac{d}{2}+1}_{2,1})}\right)^{\frac12}\right)\|\nb u\|_{L^2_t(\dot{H}^{s})}\right).
 \eeq
Using commutator estimate and product estimate \eqref{product2} in Besov space again, we get for $-\frac{d}{2}<s<\frac{d}{2}$,
 \beq\label{est_nonlin-utau}
\nn&&\sqrt{\mathrm{We}}\int_0^t\sum_{q\in\mathbb{Z}}2^{q\frac{d}{2}}\|[\ddl,u]\cdot\nb\tau\|_{L^2}+\|g_\al(\tau,\nb u)\|_{\dot{B}^{\frac{d}{2}}_{2,1}}dt'\\
\nn&\le&C\sqrt{\mathrm{We}}\int_0^t\|\nb u\|_{\dot{B}^{\frac{d}{2}}_{2,1}}\|\tau\|_{\dot{B}^{\frac{d}{2}}_{2,1}}dt'\\
\nn&\le&C\sqrt{\mathrm{We}}\left(\|\tau\|_{L^\infty_t(\dot{B}^{\frac{d}{2}}_{2,1})}\|\nb u\|_{L^1_t(\dot{B}^{\frac{d}{2}}_{2,1})}^h+\|\tau\|_{L^2_t(\dot{B}^{\frac{d}{2}}_{2,1})}\|\nb u\|_{L^2_t(\dot{B}^{\frac{d}{2}}_{2,1})}^l\right)\\
&\leq&C\sqrt{\mathrm{We}}\left(\|\tau\|_{L^\infty_t(\dot{B}^{\frac{d}{2}}_{2,1})}\|\nb u\|_{L^1_t(\dot{B}^{\frac{d}{2}}_{2,1})}^h+\left(\|\tau\|_{L^2_t(\dot{H}^s)}+\left(\|\tau\|_{L^\infty_t(\dot{B}^{\frac{d}{2}}_{2,1})}^h\|
\tau\|_{L^1_t(\dot{B}^{\frac{d}{2}}_{2,1})}^h\right)^\frac12\right)\|\nb u\|_{L^2_t(\dot{H}^s)}\right).
 \eeq
Substituting \eqref{est_conv-u-B} and \eqref{est_nonlin-utau} into \eqref{B_est} yields
  \beq\label{est_B}
\nn&&\sqrt{2\om\mathrm{Re}}\|u\|^h_{\tl{L}^\infty_t(\dot{B}^{\frac{d}{2}}_{2,1})}+\sqrt{\mathrm{We}}\|\tau\|^h_{\tl{L}^\infty_t(\dot{B}^{\frac{d}{2}}_{2,1})}\\
\nn&&\qquad\qquad\qquad+\min\left\{\sqrt{\frac{1-\om}{\mathrm{Re}}},\frac{1}{\sqrt{\mathrm{We}}}\right\}\left(\sqrt{2\om(1-\om)}\| u\|^h_{{L}^1_t(\dot{B}^{\frac{d}{2}+1}_{2,1})}+\|\tau\|^h_{{L}^1_t(\dot{B}^{\frac{d}{2}}_{2,1})}\right)\\
\nn&\le&\sqrt{2\om\mathrm{Re}}\|u_0\|^h_{\dot{B}^{\frac{d}{2}}_{2,1}}+\sqrt{\mathrm{We}}\|\tau_0\|^h_{\dot{B}^{\frac{d}{2}}_{2,1}}\\
\nn&&+C\sqrt{2\om\mathrm{Re}}\left(\|u\|_{L^\infty_t(\dot{B}^{\frac{d}{2}}_{2,1})}\|\nb u\|^h_{L^1_t(\dot{B}^{\frac{d}{2}}_{2,1})}+\left(\|u\|_{L^2_t(\dot{H}^{s+1})}+\left(\|u\|^h_{L^\infty_t(\dot{B}^{\frac{d}{2}}_{2,1})}
\|u\|^h_{L^1_t(\dot{B}^{\frac{d}{2}+1}_{2,1})}\right)^{\frac12}\right)\|\nb u\|_{L^2_t(\dot{H}^{s})}\right)\\
&&+C\sqrt{\mathrm{We}}\left(\|\tau\|_{L^\infty_t(\dot{B}^{\frac{d}{2}}_{2,1})}\|\nb u\|_{L^1_t(\dot{B}^{\frac{d}{2}}_{2,1})}^h+\left(\|\tau\|_{L^2_t(\dot{H}^s)}+\left(\|\tau\|_{L^\infty_t(\dot{B}^{\frac{d}{2}}_{2,1})}^h\|
\tau\|_{L^1_t(\dot{B}^{\frac{d}{2}}_{2,1})}^h\right)^\frac12\right)\|\nb u\|_{L^2_t(\dot{H}^s)}\right).
 \eeq
Let us define
 \beqno
E_1(t)&:=&\sqrt{\om\mathrm{Re}}\|u\|_{\tl{L}^\infty_t(\dot{H}^s)}+\sqrt{\mathrm{We}}\|\tau\|_{\tl{L}^\infty_t(\dot{H}^s)}+
\sqrt{\om(1-\om)}\|\nb u\|_{L^2_t(\dot{H}^s)}+\|\tau\|_{L^2_t(\dot{H}^s)},\\
E_2(t)&:=&\sqrt{\om\mathrm{Re}}\|u\|^h_{\tl{L}^\infty_t(\dot{B}^{\frac{d}{2}}_{2,1})}+\sqrt{\mathrm{We}}\|\tau\|^h_{\tl{L}^\infty_t(\dot{B}^{\frac{d}{2}}_{2,1})}+
\sqrt{\om(1-\om)}\|\nb u\|^h_{L^1_t(\dot{B}^{\frac{d}{2}}_{2,1})}+\|\tau\|^h_{L^1_t(\dot{B}^{\frac{d}{2}}_{2,1})},\\
 \eeqno
 \beno
E_1(0):=\sqrt{\om\mathrm{Re}}\|u_0\|_{\dot{H}^s}+\sqrt{\mathrm{We}}\|\tau_0\|_{\dot{H}^s}, \qquad E_2(0):=\sqrt{\om\mathrm{Re}}\|u_0\|^h_{\dot{B}^{\frac{d}{2}}_{2,1}}+\sqrt{\mathrm{We}}\|\tau_0\|^h_{\dot{B}^{\frac{d}{2}}_{2,1}},
 \eeno
and
$$
E(t):=E_1(t)+E_2(t), \qquad E(0):=E_1(0)+E_2(0).
$$
Moreover, we denote
 \beqno
\kappa_1&:=&\max\left\{\left(\frac{\mathrm{We}}{\om(1-\om)}\right)^{\frac14}, \frac{1}{\left(\om(1-\om)\right)^{\frac14}}, \frac{\left(\om\mathrm{Re}\right)^\frac14}{\sqrt{\om(1-\om)}}, \frac{\left(\om\mathrm{Re}\right)^{\frac18}}{\left(\om(1-\om)\right)^{\frac38}}\right\},\\
\kappa_2&:=&\max\left\{1, \sqrt{\frac{\mathrm{Re}}{1-\om}}, \sqrt{\mathrm{We}}\right\},\\
\kappa_3&:=&\max\left\{\frac{1}{\sqrt{\om(1-\om)}}, \frac{\sqrt{\mathrm{Re}}}{\sqrt{\om}(1-\om)}, \frac{(\mathrm{Re})^\frac14}{\om^{\frac12}(1-\om)^{\frac34}}, \frac{(\mathrm{We})^\frac14}{\sqrt{\om(1-\om)}}\right\}.
 \eeqno
Then \eqref{est_H^s} and \eqref{est_B} read as follows:
 \be
E_1(t)\leq E_1(0)+C\kappa_1E(t)^{\frac32}, \quad \hbox{for} \quad -\frac{d}{2}<s<\frac{d}{2}-1,
 \ee
and
 \be
E_2(t)\leq \kappa_2E_2(0)+C\kappa_2\kappa_3E(t)^{2}, \quad \hbox{for} \quad -\frac{d}{2}<s<\frac{d}{2}-1,
 \ee
Consequently,
 \be\label{eq3.19}
E(t)\leq \kappa_2E(0)+C\left(\kappa_1E(t)^{\frac32}+\kappa_2\kappa_3E(t)^{2}\right), \quad \hbox{for} \quad -\frac{d}{2}<s<\frac{d}{2}-1.
 \ee
By using standard continuity method, we infer from  \eqref{eq3.19} that
 \be\label{eq3.20}
E(t)\leq 2\kappa_2E(0),
 \ee
provided $E(0)$ is small enough.  Then the existence part of Theorem \ref{thm-global} follows immediately. $\quad\quad\,\Box$
\section{Uniqueness}\label{Sec-U}
\noindent Let $(u_1,\tau_1)$ and $(u_2,\tau_2)$ be the solution to the system \eqref{IOBdimensionless} with the same initial data obtained in Section \ref{Sec-E}. Denote $(w, \sig):=(u_1-u_2, \tau_1-\tau_2)$, and $p=\Pi_1-\Pi_2$. Then it is easy to verify that $(w,\sig)$ satisfies the following system:
\beq\label{difference}
\begin{cases}
\Rey\pr_tw-(1-\om)\Dl w+\nb p=\dv \sig-\Rey(w\cdot\nb u_1+u_2\cdot \nb w),\\
\We(\pr_t\sig+u_1\cdot\nb\sig)+\sig=2\om D(w)-\We w\cdot\nb\tau_2-\We g_\al(\sig,\nb u_1)-\We g_{\al}(\tau_2,\nb w).
\end{cases}
\eeq
Applying the localized operator $\ddl$  to system \eqref{difference} yields
\beq\label{l-difference}
\begin{cases}
\Rey\pr_t\ddl w-(1-\om)\Dl \ddl w+\nb \ddl p=\dv \ddl\sig-\Rey\ddl(w\cdot\nb u_1+u_2\cdot \nb w),\\
\We(\pr_t\ddl\sig+u_1\cdot\nb\ddl\sig)+\ddl\sig=2\om D(\ddl w)-\We[\ddl,u_1]\cdot\nb\sig-\We \ddl(w\cdot\nb\tau_2)\\
\qquad\qquad\qquad\qquad\qquad\qquad\qquad-\We \ddl g_\al(\sig,\nb u_1)-\We\ddl g_{\al}(\tau_2,\nb w).
\end{cases}
\eeq
Using the cancelation relation $(\dv\ddl\sig|\ddl w)+(D(\ddl w)|\ddl \sig)=0$, similar to \eqref{est_loc}, we arrive at
 \beq\label{est_uni1}
\nn&&\frac12\frac{d}{dt}\left(2\om\mathrm{Re}\|\ddl w\|_{L^2}^2+\mathrm{We}\|\ddl\sig\|_{L^2}^2\right)+2\om(1-\om)\|\nb\ddl w\|_{L^2}^2+\|\ddl\sig\|_{L^2}^2\\
\nn&\le&2\om\mathrm{Re}\|\ddl(w\cdot\nb u_1+u_2\cdot\nb w)\|_{L^2}\|\ddl w\|_{L^2}+\mathrm{We}\left(\|[\ddl,u_1]\cdot\nb\sig\|_{L^2}+\|\ddl(w\cdot\nb\tau_2)\|_{L^2}\right)\|\ddl \sig\|_{L^2}\\
&&+\We\left(\|\ddl g_\al(\sig,\nb u_1)\|_{L^2}+\|\ddl g_\al(\tau_2,\nb w)\|_{L^2}\right)\|\ddl \sig\|_{L^2}.
 \eeq
Integrating  \eqref{est_uni1} w. r. t. time $t$, multiplying the resulting inequality by $2^{2qs}$, and then taking sum w. r. t. $q$ over $\Z$, using H\"{o}lder's inequality, we are led to
 \beq\label{est_uni2}
\nn&&\om\mathrm{Re}\| w(t)\|_{\dot{H}^s}^2+\fr{\mathrm{We}}{2}\|\sig(t)\|_{\dot{H}^s}^2+2\om(1-\om)\|\nb w\|_{L^2_t(\dot{H}^s)}^2+\|\sig\|_{L^2_t(\dot{H}^s)}^2\\
\nn&\le&2\om\mathrm{Re}\int_0^t\|w\cdot\nb u_1+u_2\cdot\nb w\|_{\dot{H}^s}\| w\|_{\dot{H}^s}dt'+\mathrm{We}\int_0^t\left(\sum_{q\in\Z}2^{2qs}\|[\ddl,u_1]\cdot\nb\sig\|_{L^2}^2\right)^{\fr12}\|\sig\|_{\dot{H}^s}dt'\\
&&+\We\int_0^t\left(\|w\cdot\nb\tau_2\|_{\dot{H}^s}+\| g_\al(\sig,\nb u_1)\|_{\dot{H}^s}+\| g_\al(\tau_2,\nb w)\|_{\dot{H}^s}\right)\|\sig\|_{\dot{H}^s}dt'.
 \eeq
By virtue of the product estimates \eqref{product1}, \eqref{product1.5} and commutator estimate in Besov spaces, we have
\begin{gather}
\nn\|w\cdot\nb u_1+u_2\cdot\nb w\|_{\dot{H}^s}\lesssim \|\nb u_1\|_{\dot{B}^{\fr{d}{2}}_{2,1}}\|w\|_{\dot{H}^s}+\|u_2\|_{\dot{B}^{\fr{d}{2}}_{2,1}}\|\nb w\|_{\dot{H}^s},\quad \mathrm{for}\quad -\fr{d}{2}<s<\fr{d}{2},\\
\nn\left(\sum_{q\in\Z}2^{2qs}\|[\ddl,u_1]\cdot\nb\sig\|_{L^2}^2\right)^{\fr12}\lesssim \|\nb u_1\|_{\dot{B}^{\fr{d}{2}}_{2,1}}\|\sig\|_{\dot{H}^s}, \quad\mathrm{for}\quad  -\fr{d}{2}-1<s<\fr{d}{2},\\
\nn\|w\cdot\nb\tau_2\|_{\dot{H}^s}\lesssim \|w\|_{\dot{H}^{s+1}}\|\nb\tau_2\|_{\dot{B}^{\fr{d}{2}-1}_{2,\infty}}\lesssim\|\nb w\|_{\dot{H}^{s}}\|\tau_2\|_{\dot{B}^{\fr{d}{2}}_{2,1}}, \quad\mathrm{for}\quad-\fr{d}{2}<s<\fr{d}{2}-1,\\
\nn\| g_\al(\sig,\nb u_1)\|_{\dot{H}^s}+\| g_\al(\tau_2,\nb w)\|_{\dot{H}^s}\lesssim\|\nb u_1\|_{\dot{B}^{\fr{d}{2}}_{2,1}}\|\sig\|_{\dot{H}^s}+\|\tau_2\|_{\dot{B}^{\fr{d}{2}}_{2,1}}\|\nb w\|_{\dot{H}^s},\quad \mathrm{for}\quad -\fr{d}{2}<s<\fr{d}{2},
\end{gather}
Substituting these estimates into \eqref{est_uni2} yields that, for $-\fr{d}{2}<s<\fr{d}{2}-1$, there holds
 \beq\label{est_uni3}
\nn&&\om\mathrm{Re}\| w(t)\|_{\dot{H}^s}^2+\mathrm{We}\|\sig(t)\|_{\dot{H}^s}^2+2\om(1-\om)\|\nb w\|_{L^2_t(\dot{H}^s)}^2+\|\sig\|_{L^2_t(\dot{H}^s)}^2\\
\nn&\le&C\int_0^t\|\nb u_1\|_{\dot{B}^{\fr{d}{2}}_{2,1}}\left(\om\Rey\|w\|_{\dot{H}^s}^2+\We\|\sig\|_{\dot{H}^s}^2\right)dt'\\
&&+C\int_0^t\left(\om\Rey\|u_2\|_{\dot{B}^{\fr{d}{2}}_{2,1}}\|w\|_{\dot{H}^s}+\We\|\tau_2\|_{\dot{B}^{\fr{d}{2}}_{2,1}}\|\sig\|_{\dot{H}^s}\right)\|\nb w\|_{\dot{H}^s}dt'.
 \eeq
Noting that by Cauchy's inequality, there exists a positive constant $C$ depending on $\Rey$, $\We$ and $\om$, such that
\beqno
&&\int_0^t\left(\om\Rey\|u_2\|_{\dot{B}^{\fr{d}{2}}_{2,1}}\|w\|_{\dot{H}^s}+\We\|\tau_2\|_{\dot{B}^{\fr{d}{2}}_{2,1}}\|\sig\|_{\dot{H}^s}\right)\|\nb w\|_{\dot{H}^s}dt'\\
&\le&\om(1-\om)\|\nb w\|_{L^2_t(\dot{H}^s)}^2+C\int_0^t\left(\om\Rey\|u_2\|_{\dot{B}^{\fr{d}{2}}_{2,1}}^2+
\We\|\tau_2\|_{\dot{B}^{\fr{d}{2}}_{2,1}}^2\right)\left(\om\Rey\|w\|_{\dot{H}^s}^2+\We\|\sig\|_{\dot{H}^s}^2\right)dt',
\eeqno
combining this inequality with \eqref{est_uni3}, we obtain
 \beq\label{est_uni4}
\nn&&\om\mathrm{Re}\| w(t)\|_{\dot{H}^s}^2+\mathrm{We}\|\sig(t)\|_{\dot{H}^s}^2+\om(1-\om)\|\nb w\|_{L^2_t(\dot{H}^s)}^2+\|\sig\|_{L^2_t(\dot{H}^s)}^2\\
\nn&\le&C\int_0^t\left(\|\nb u_1\|_{\dot{B}^{\fr{d}{2}}_{2,1}}+\om\Rey\|u_2\|_{\dot{B}^{\fr{d}{2}}_{2,1}}^2+
\We\|\tau_2\|_{\dot{B}^{\fr{d}{2}}_{2,1}}^2\right)\left(\om\Rey\|w\|_{\dot{H}^s}^2+\We\|\sig\|_{\dot{H}^s}^2\right)dt'.
 \eeq
Thanks to  the embedding in low frequency,  we infer from \eqref{eq3.20} that
\beqno
&&\int_0^t\left(\|\nb u_1\|_{\dot{B}^{\fr{d}{2}}_{2,1}}+\om\Rey\|u_2\|_{\dot{B}^{\fr{d}{2}}_{2,1}}^2+
\We\|\tau_2\|_{\dot{B}^{\fr{d}{2}}_{2,1}}^2\right)dt'\\
&\lesssim&\|\nb u_1\|^h_{L^1_t(\dot{B}^{\fr{d}{2}}_{2,1})}+t^{\fr12}\|\nb u_1\|_{L^2_t(\dot{H}^{s})}+\om\Rey t\left(\|u_2\|^h_{\tl{L}^{\infty}_t(\dot{B}^{\fr{d}{2}}_{2,1})}+\|u_2\|_{\tl{L}^{\infty}_t(\dot{H}^{s})}\right)^2\\
&&+\We t\left(\|\tau_2\|^h_{\tl{L}^{\infty}_t(\dot{B}^{\fr{d}{2}}_{2,1})}+\|\tau_2\|_{\tl{L}^{\infty}_t(\dot{H}^{s})}\right)^2\\
&\le&\infty,
\eeqno
for any $t>0$. Then the uniqueness follows from Glonwall's inequality immediately. This completes the proof of Theorem \ref{thm-global}.$\Box$
 \bigbreak\noindent{\bf Acknowledgment.}
This work is supported by China Postdoctoral Science Foundation funded project 2014M552065, and National Natural Science Foundation of China 11401237.

\end{document}